\documentclass[letterpaper, 10 pt, conference]{ieeeconf}  

\IEEEoverridecommandlockouts                              
\usepackage{graphicx} 
\usepackage{amsmath} 
\usepackage{amssymb}  
\usepackage{color,soul}
\usepackage[dvipsnames]{xcolor}
\usepackage{algorithm}
\usepackage{mathtools}
\usepackage{algpseudocode}
\algtext*{EndFor}
\algtext*{EndIf}

\usepackage{subcaption} 
\usepackage{hyperref}
\usepackage{enumerate}
\usepackage{cite}

\usepackage{tikz}
\usetikzlibrary{calc, shapes, backgrounds}
\usepackage{pgfpages}
\DeclareMathAlphabet{\mathcal}{OMS}{cmsy}{m}{n}
\usepackage{pgf}
\usetikzlibrary{arrows,automata,calc,positioning}

\usepackage{xcolor}
\definecolor{myLightBlue}{rgb}{0.357,0.627,0.796}
\definecolor{myMidBlue}{rgb}{0,0.4,0.663}
\definecolor{matlabBlue}{rgb}{0,0.447,0.741}
\definecolor{myRed}{rgb}{0.918,0.133,0.047}
\definecolor{myGrey}{rgb}{0.75,0.75,0.75}
\definecolor{myLightGrey}{rgb}{0.9,0.9,0.9}
\definecolor{myOrange}{rgb}{0.85,0.325,0.098}
\definecolor{myGreen}{rgb}{0.466,0.674,0.188}
\definecolor{myPurple}{rgb}{0.5,0.3,1}
\definecolor{mygreen}{RGB}{204,255,153}
\definecolor{myblue}{RGB}{153,204,255}
\definecolor{myorange}{RGB}{255,204,153}
\definecolor{mylightgray}{RGB}{225,225,225}

\newtheorem{definition}{\bf Definition}
\newtheorem{proposition}{\bf Proposition}[section]
\newtheorem{lemma}{\bf Lemma}[section]

\newtheorem{theorem}{\bf Theorem}[section]

\renewenvironment{proof}{{\bfseries Proof:}}{\cqfd}

\providecommand{\com[2]}{\begin{tt}[#1: #2]\end{tt}}

\providecommand{\prt}[1]{\ensuremath{\left( #1 \right)}}
\newcommand{\cqfd}{\hfill \rule{2mm}{2mm}\medbreak\indent}

\DeclareMathOperator{\newrank}{rank} 
\newcommand{\rank}{\newrank\,}

\DeclareMathOperator{\newsgn}{sgn} 
\newcommand{\sgn}{\newsgn}

\newcommand{\toadd}[1]{}

\newcommand{\Gleft}{G_B}
\newcommand{\Gright}{G_C}

\newcommand{\Bleft}{B^\dag}
\newcommand{\Cright}{C^\dag}

\title{\LARGE \bf
Combinatorial Characterization for Global Identifiability\\
of Separable Networks with Partial Excitation and Measurement
}

\author{Antoine Legat and Julien M. Hendrickx
\thanks{A. Legat and J. M. Hendrickx are with ICTEAM Institute, UCLouvain, Belgium.
Work supported
by the ``RevealFlight'' ARC at UCLouvain, by the Incentive Grant
for Scientific Research (MIS) ``Learning from Pairwise Data'' and by the KORNET project from F.R.S.-FNRS {\tt\small antoine.legat@uclouvain.be, julien.hendrickx@uclouvain.be.}}%
}

\begin{document}

\maketitle
\thispagestyle{empty}
\pagestyle{empty}

\begin{abstract}
This work focuses on the generic identifiability of dynamical networks with partial excitation and measurement: a set of nodes are interconnected by transfer functions according to a known topology, some nodes are excited, some are measured, and only a part of the transfer functions are known. Our goal is to determine whether the unknown transfer functions can be generically recovered based on the input-output data collected from the excited and measured nodes.

We introduce the notion of separable networks, for which global and so-called local identifiability are equivalent.
A novel approach yields a necessary and sufficient combinatorial characterization for local identifiability for such graphs, in terms of existence of paths and conditions on their parity.
Furthermore, this yields a necessary condition not only for separable networks, but for networks of any topology.
\end{abstract}

\smallskip

%
%
%
\section{INTRODUCTION}
\medskip

This paper addresses the identifiability of dynamical networks in which node signals are connected by causal linear time-invariant transfer functions, and can be excited and/or measured. Such networks can be modeled as directed graphs where each edge carries a transfer function, and known excitations and measurements are applied at certain nodes.

\smallskip

\subsection{Problem Statement}


We consider the identifiability of a network matrix $G(q)$, where the network is made up of $n$ node signals stacked in the vector $w(t) = [w_1(t) \quad w_2(t) \ \cdots \ w_n(t)]^\top$, known external excitation signals $r(t)$, measured node signals $y(t)$ and unmeasured noise $v(t)$ related to each other by:
\begin{align} \label{eq:networkModel} \begin{split}
	w(t) &= G(q)w(t) + Br(t) + v(t)\\
	y(t) &= Cw(t),
\end{split} \end{align}
where matrices $B$ and $C$ are binary selections indicating respectively the $n_B$ excited and $n_C$ measured nodes, forming sets $\mathcal{B}$ and $\mathcal{C}$ respectively. Matrix $B$ is full column rank and each column contains one $1$ and $n-1$ zeros. Matrix $C$ is full row rank and each row contains one $1$ and $n-1$ zeros.

\medskip

The nonzero entries of the transfer matrix $G(q)$ define the network topology: $G_{ij}(q)$ is the transfer function from node $j$ to node $i$.
It is represented by an edge $(j,i) \in \mathcal E$ in the graph, where $\mathcal E$ is the set of all edges, each corresponding to a nonzero entry of $G(q)$.
Some of those transfer functions are known and collected in the matrix $G^\bullet(q)$, and the others unknown, collected in matrix $G^\circ(q)$, such that $G(q) = G^\bullet(q) + G^\circ(q)$. The known edges (i.e. the edges corresponding to known transfer functions) are collected in set $\mathcal E^\bullet$, the unknown ones in $\mathcal E^\circ$, and they form a partition of the set of all edges $\mathcal E = \{ \mathcal E^\bullet, \mathcal E^\circ\}$. We denote the number of unknown transfer functions by $m_\circ \triangleq |\mathcal E^\circ|$. \toadd{Add example to help with notations?}

\medskip
		
We assume that \emph{the input-output relations between the excitations $r$ and measurements $y$ have been identified}, and that the network topology is known.
From this knowledge, we aim at recovering the unknown transfer functions $G^\circ(q)$.

\smallskip

\subsection{State of the Art}

\smallskip

The model \eqref{eq:networkModel} has recently been the object of a significant research effort.
Network identifiability was first introduced in \cite{weerts2015identifiability}, in case the whole network is to be recovered.
Conditions for the identification of a single transfer function are derived in \cite{weerts2018single, gevers2018practical}.
Studying the influence of rank-reduced or correlated noise under certain assumptions yields less conservative identifiability conditions \cite{weerts2018prediction, gevers2018identifiability, van2019local, ramaswamy2020local}.
	
\medskip

It turns out that under some assumptions, identifiability of the network, i.e. the ability to recover a transfer function or the whole network from the input-output relation, is a generic notion: Either \emph{almost all} transfer matrices corresponding to a given network structure are identifiable, in which case the structure is called \emph{generically identifiable}, or none of them are. 
A number of works study generic identifiability when all nodes are excited or all nodes are measured, i.e. when $B$ or $C=I$ \cite{bazanella2017identifiability, weerts2018identifiability}. Considering the graph of the network, path-based conditions on the allocation of measurements (resp. excitations) in the case of full excitation (resp. measurement) are derived in \cite{hendrickx2018identifiability} (resp. \cite{shi2020excitation}). Reformulating these conditions by means of disjoint trees in the graph, the authors of \cite{cheng2019allocation, cheng2021allocation, dreef2022excitation}
develop scalable algorithms to allocate excitations/measurements in case of full measurement/excitation. In case of full measurement, \cite{shi2020generic} derives path-based conditions for the generic identifiability of a subset of transfer functions, with noise. 
	
\medskip

While the conditions in \cite{hendrickx2018identifiability, shi2020excitation} apply to generic identifiability i.e. for almost all transfer functions, \cite{van2018topological} extends these results to the stronger requirement of identifiability \emph{for all} (nonzero) transfer matrices corresponding to a given structure, and \cite{van2019necessary} provides conditions for the outgoing edges of a node, and the whole network under the same conditions.
	
\medskip
	
As mentioned, the common assumption in all these works is that is that either all nodes are excited, or they are all measured. In \cite{bazanella2019network}, this assumption is relaxed and generic identifiability with partial excitation and measurement is addressed for particular network topologies.
For acyclic networks, \cite{cheng2022identifiability} gives necessary conditions, introduces the transpose network and shows that identifiability of the transpose network and its original network is equivalent.
	
\medskip

For an arbitrary topology, \cite{shi2021single} provides for noise exploitation an elegant reformulation as an equivalent network model, where noise is cast into excitation signals.
\cite{cheng2021necessary, mapurunga2022identifiability, mapurunga2022excitation} provide necessary conditions for network identifiability, but do not handle a priori known/fixed transfer functions, which could lead to less conservative conditions.

\medskip

In the general case of arbitrary topology, partial excitation and measurement, we introduced in \cite{legat2020local} the notion of local identifiability, i.e. only on a neighborhood of $G(q)$. Local identifiability is a generic property, necessary for generic identifiability and no counterexample to sufficiency is known to the authors, i.e. no network which is locally identifiable but not globally identifiable. We derived algebraic necessary and sufficient conditions for generic local identifiability for both the whole network, and a single transfer function.

\medskip

The algebraic conditions of \cite{legat2020local} allow rapidly testing local identifiability for any given network,
but finding a graph-theoretical characterization, akin to what was done in the full excitation case \cite{hendrickx2018identifiability}, remains an open question, even though it is known that the property depends solely on the graph.
Such characterization would in particular pave the way for optimizing the selection of nodes to be excited and measured, alike the work in \cite{shi2020excitation} in the full measurement case.

\medskip

In this line of work, we introduced in \cite{legat2021path} \emph{decoupled} identifiability, necessary for local identifiability, see Section \ref{sec:identif}. 
We extended the algebraic characterization of \cite{legat2020local} when some transfer functions are known/fixed a priori, and developed it in terms of closed-loop transfer matrices $T(G)=(I-G)^{-1}$, which led to some necessary and some sufficient path-based conditions for decoupled identifiability.


An approach different to all that precedes is to that network dynamics are known, and aim at identifying the topology from input/output data. This problem is referred to as topology identification, and is addressed in e.g. \cite{van2019topology_heterogeneous, van2019topology_lyapunov}. \cite{kivits2022local} studies diffusively coupled linear networks, which can be represented by undirected graphs.

\medskip

In this paper, we introduce \emph{separable networks}, which are a generalization of the decoupled version of the network, allowing for different topologies in excited and measured subgraphs, see Section \ref{sec:separable} and Fig. \ref{fig:separable_ex}. Thanks to their particular structure, global and local identifiability are equivalent on those networks, thus global identifiability can be studied with the algebraic tools of \cite{legat2020local}.

\smallskip

We obtain a necessary and sufficient \emph{combinatorial characterization} of global identifiability, in terms of existence of paths and conditions on their parity.
Since the decoupled network is a particular case of separable network, this condition can be formulated on the decoupled network of networks of any topology, \emph{not only separable networks}. And generic decoupled identifiability is necessary for generic global identifiability \cite{legat2021path}, hence the necessary condition applies to a network of any topology.



\smallskip

\subsection{Framework} \label{subsec:framework}

\smallskip
\noindent \textbf{Assumptions: }
	We consider model \eqref{eq:networkModel}. Consistently with previous work (e.g.\cite{hendrickx2018identifiability, van2019necessary, weerts2018identifiability}), we assume:
	\begin{enumerate}[1)]
		\item The network is well-posed and stable, that is $(I-G(q))^{-1}$ is proper and stable.
		\item All the entries of $G(q)$ are proper transfer functions. 
        \item The network is stable in the following sense : $|\lambda_i| < 1$ for each eigenvalue $\lambda_i$ of $G$. \label{assumption:eigenvalues_G}
	\end{enumerate}

\medskip

Throughout the paper, we develop our results without exploiting noise signals.
However, under some mild assumptions, noise signals $v(t)$ can play the same role as excitation signals, as \cite{shi2021single}: the network is reformulated as an equivalent network model, where noise is cast into excitations.

\medskip

We suppose that there are exactly $n_B n_C$ unknown transfer functions, i.e. as many as the number of (excitation - measurement) pairs, and we address the identifiability of the \emph{whole} network, i.e. all unknown transfer functions.
If there are more unknown transfer functions, then it is not identifiable since there are more unknowns than (input, output) data.


\medskip

\noindent\textbf{Genericity:}
We will focus on \emph{generic} properties: we say that a property is generic if it either holds (i) for \emph{almost all} variables, i.e for all variables except possibly those lying on a lower-dimensional set\cite{davison1977connectability, dion2003generic} (i.e. a set of dimension lower than $m_\circ$, the number of unknown transfer functions), or (ii) for no variable.
For example, take a polynomial $p$. The \emph{nonzeroness} of $p(x)$ is a generic property of $x$: either (i) $p(x) \neq 0$ for all $x$ except its roots, or (ii) $p$ is the zero polynomial, which returns zero for all $x$.

\medskip

Consistently with \cite{legat2020local, legat2021path}, we consider a single frequency:
instead of working with transfer functions $G_{ij}(q)$ and their transfer matrix $G(q)$, we work with scalar values $G_{ij}$ and their matrix $G \in \mathbb C^{n \times n}$.
Conceptually, our generic results directly extend to the transfer function case: if one can recover a $G_{ij}(z)$ at a given frequency $z$ for almost all $G$ consistent with a network, then one can also recover it at all other frequencies.
In the remainder, for the sake of simplicity, we work in the scalar setup, hence omit the $(q)$.

\bigskip

\noindent The framework of this paper is summarized below:

\begin{itemize}
    \item Partial excitation and measurement
    \item Allows for the presence of known transfer functions next to the unknown ones
    \item No use of noise, scalar setup
	\item \emph{Global} identifiability on \emph{separable} networks
	\begin{itemize}
		\item[$\Rightarrow$] Applies on decoupled identifiability of all networks
	\end{itemize}
	\item \emph{Generic} identifiability of the \emph{whole} network
	\item We assume $n_B n_C = m_\circ$, i.e. as much data as unknowns
\end{itemize}

\medskip


\medskip

\newpage

\section{IDENTIFIABILITY} \label{sec:identif}
\medskip

We remind different notions of identifiability relevant to our work \cite{weerts2015identifiability, hendrickx2018identifiability, legat2020local, legat2021path}. 
In order to lighten notations, we denote $T(G) = (I-G)^{-1}$. 

\medskip

\begin{definition} \label{def:global_identif}
	A network is \emph{globally identifiable} at $G$ from excitations $\mathcal{B}$ and measurements $\mathcal{C}$ if, for all network matrix $\tilde{G}$ with same zero and known entries as $G$, there holds
	\begin{align} \label{eq:def_global_identif}
		C \, T(\tilde{G}) \, B = C \, T(G) \, B \Rightarrow \tilde{G}^\circ= G^\circ,
	\end{align}
    where matrices $B$ and $C$ are binary selections indicating respectively the excited and measured nodes, forming sets $\mathcal{B}$ and $\mathcal{C}$ respectively. $C \, T(G) \, B$ is the global transfer matrix between input and output.
	The network is \emph{generically globally identifiable} if it is identifiable at \emph{almost all} $G$.
\end{definition}

\bigskip

This definition extends \cite{hendrickx2018identifiability} to the case where some transfer functions are known ($G^\bullet$), and some are not ($G^\circ$), as in \cite{weerts2018single}.
We call \emph{global identifiability} this standard notion of identifiability, to avoid confusion with \emph{local} identifiability \cite{legat2020local}, which corresponds to identifiability provided that $\tilde G$ is sufficiently close to $G$. Local identifiability is necessary for global identifiability, and no counter-example to sufficiency is known.

\medskip

\begin{definition} \label{def:local_identif}
	A network is \emph{locally identifiable} at $G$ from excitations $\mathcal{B}$ and measurements $\mathcal{C}$ if there exists $\epsilon > 0$ such that for any $\tilde{G}$ with same zero and known entries as $G$ satisfying $||\tilde{G}-G||<\epsilon$, there holds
	\begin{align} \label{eq:def_local_identif}
		C \, T(\tilde{G}) \, B = C \, T(G) \, B \Rightarrow \tilde{G}^\circ= G^\circ.
	\end{align}
	The network is \emph{generically locally identifiable} if it is locally identifiable at \emph{almost all} $G$.
\end{definition}

\bigskip

Local identifiability can be characterized algebraically based on the matrix $K$:
\begin{align} \label{eq:def_K}
	K(G) \triangleq \prt{B^\top T^\top(G) \otimes C \, T(G)} I_{G^\circ},
\end{align}
where we remind that $T(G) = (I-G)^{-1}$, symbol $\otimes$ denotes the Kronecker product and the matrix $I_{G^\circ} \in \{0,1\}^{n^2 \times m_\circ}$ selects only the columns of the preceding $n_B n_C \times n^2$ matrix corresponding to unknown modules \cite{legat2021path}.

\medskip

\begin{theorem} \label{thm:network_local_identif}
        (Corollary 4.1 in \cite{legat2020local})\\
	    Exactly one of the two following holds:
	    \begin{enumerate}[(i)]
		    \item $\rank K = m_\circ$ for almost all $G$ and $G^\circ$ is locally identifiable at almost all $G$;
		    \item $\rank K < m_\circ$ for all $G$ and $G^\circ$ is locally \emph{non-}identifiable at all $G$, therefore globally non-identifiable at all $G$.
		\end{enumerate}
		Moreover, $\rank K = m_\circ$ is equivalent to the following implication holding for all $\Delta$ with same zero entries as $G^\circ$:
		\begin{align} \label{eq:CTDeltaTB_Delta_net}
			C \, T(G) \, \Delta \, T(G) \, B = 0 \Rightarrow \Delta = 0.
		\end{align}
\end{theorem}

\bigskip

Equation \eqref{eq:CTDeltaTB_Delta_net} suggests the definition of a new notion, where the two matrices $T(G)$ of \eqref{eq:CTDeltaTB_Delta_net} do not need to have the same parameters anymore: this notion is \emph{decoupled} identifiability.

\medskip

\begin{definition} \label{def:decoupled_identif}
    A network is \emph{decoupled-identifiable} at $(G,$ $G')$, with $G$ and $G'$ sharing the same zero entries, if for all $\Delta$ with same zero entries as $G^\circ$, there holds:
    \begin{align} \label{eq:CTDeltaT'B_bili}
		C  \, T(G) \,  \Delta \,  T(G')  \, B = 0 \Rightarrow \Delta = 0.
	\end{align}
\end{definition}

\medskip


Decoupled identifiability, initially introduced for purely algebraic reasons, can be interpreted in terms of the identifiability of a larger network: the \emph{decoupled network.}

\medskip

\begin{definition} \label{def:decoupled_network}
    Consider a network of $n$ nodes with excitation matrix $B$, measurement matrix $C$ and network matrix $G = G^\bullet + G^\circ$, where $G^\bullet$ collects the known modules and $G^\circ$ collects the unknown modules. Its \emph{decoupled network} is composed of $2n$ nodes: $\{1, \dots, n, 1', \dots, n'\}$. Its network matrix is defined by
    \begin{align*}
        \hat G (G,G') \triangleq
        \begin{bmatrix}
            G & G^\circ \\
            0 & G'
        \end{bmatrix},
    \end{align*}
    where $G'$ has the same zero entries as $G$. Transfer matrices $G$ and $G'$ are considered as known, while $G^\circ$ contains the unknown modules.
    Excitations are applied on the first subgraph ($G'$), and measurements on the second one ($G$), i.e. its excitation and measurement matrices are
    \begin{align*}
        \hat B \triangleq
        \begin{bmatrix}
            0 & 0 \\
            0 & B
        \end{bmatrix},
        && \hat C \triangleq
        \begin{bmatrix}
            C & 0 \\
            0 & 0
        \end{bmatrix}.
    \end{align*}
\end{definition}

\medskip


The proposition below relates the notion of decoupled identifiability with the decoupled network we just introduced.

\medskip

\begin{proposition} \label{prop:equivalence_bili_decoupled}
    (Proposition 3.3 in \cite{legat2021path}) 
    The network $G$ is generically decoupled-identifiable if and only if its decoupled network $\hat G$ is generically globally identifiable. 
\end{proposition}

\medskip

Generic decoupled identifiability is necessary for generic local identifiability (which is itself necessary for generic global identifiability) \cite{legat2021path}. No counterexample to sufficiency is known to the authors, despite extensive systematic numerical tests, available in \cite{matlab}.

\medskip

\begin{proposition} \label{prop:bili_necessary}
    (Proposition 3.2 in \cite{legat2021path})
    If a network is generically locally identifiable, then it is generically decoupled-identifiable.
\end{proposition}



\medskip

\section{SEPARABLE NETWORKS} \label{sec:separable}
\smallskip

We now introduce separable networks, which are a generalization of decoupled networks, where the excited and measured subgraphs are not required to have the same topology anymore.
A \emph{separable network} is a network for which excitations and measurement can be isolated in two distinct subgraphs, with no known transfer function linking the two subgraphs. Between the two subgraphs lie all the unknown transfer functions, going from the excited subgraph to the measured one.
A formal definition in terms of network matrices is given below, and an example is given in Fig. \ref{fig:separable_ex}.

\begin{figure}[ht]
\centering
		\begin{tikzpicture}[->,>=stealth',shorten >=1pt,auto,node distance=2.9cm,line width=0.4mm, nodes = {draw, circle}]
		  \node (al)						{};
		  \node (ar) [right = 1 cm of al]			{};
		  \node (bl) [below = 0.75 cm of al]			{};
		  \node (br) [right = 1 cm of bl]			{};
		  \node (cl) [below = 0.75 cm of bl]			{};
		  \node (cr) [right = 1 cm of cl]			{};
		  \node (dl) [below = 0.75 cm of cl]			{};
		  \node (dr) [right = 1 cm of dl]			{};
		  \node (ml) [left = 1 cm of bl]			{};
		  \node (mr) [right = 1 cm of ar]			{};
        \node (nl) [below = 1 cm of ml] {};
        \node (nr) [below = 1.75 cm of mr] {};
		  \node (N1) [left of=al] 	  {};
		  \node (N2) [left of=dl] 	  {};
		  \node (N3) [right of=ar] 	  {};
		  \node (N4) [right of=dr] 	  {};
		  \path   (N1) edge  (ml)
				(N2) edge  (ml)
				(ml) edge  (bl)
				(ml) edge  (cl)
				(br) edge  (mr)
				(ar) edge  (mr)
				(mr) edge  (N3)
				(mr) edge  (N4)
		  	  	(al) edge [dashed, color = myOrange]   (ar)
		  		(bl) edge [dashed, color = myOrange]   (br)
				(cl) edge [dashed, color = myOrange]   (cr)
				(dl) edge [dashed, color = myOrange]   (dr)
		 	   	(N1) edge  (al)
				(N2) edge  (dl)
				(cr) edge  (N3)
				(dr) edge  (N4)
                (N2) edge (nl)
                (ml) edge (nl)
                (nl) edge (dl)
                (cl) edge (nl)
                (N4) edge (nr)
                (dr) edge (nr)
                (nr) edge (cr)
                (nr) edge (mr)
                ;
        \draw[myMidBlue,dotted,line width=0.5mm] ($(N3.north west)+(-0.5,0.25)$)  rectangle ($(N4.south east)+(0.5,-0.25)$);
   	    \node (C) [below = 1.25 cm of N3, draw=none] {\textcolor{myMidBlue}{$\mathcal{C}$}};
   	    \draw[myMidBlue,dotted,line width=0.5mm] ($(N1.north west)+(-0.5,0.25)$)  rectangle ($(N2.south east)+(0.5,-0.25)$);
   	    \node (B) [below = 1.25 cm of N1, draw=none] {\textcolor{myMidBlue}{$\mathcal{B}$}};
        \node (G_B) [above = 10 mm of ml, draw=none] {$\Gleft$};
        \node (G_C) [above = -1 mm of mr, draw=none] {$\Gright$};
        \node (bmid) [right = 3 mm of bl, draw=none]			{};
        \node (Gcirc) [above = 10 mm of bmid, draw=none] {\textcolor{myOrange}{$G^\odot$}};
		\end{tikzpicture}
\caption{An example of separable network: the excitations $\mathcal B$ are isolated in one subgraph, and the measurements $\mathcal C$ in another subgraph. The unknown transfer functions, in dashed orange, link the excited subgraph to the measured one.}
\label{fig:separable_ex}
\end{figure}
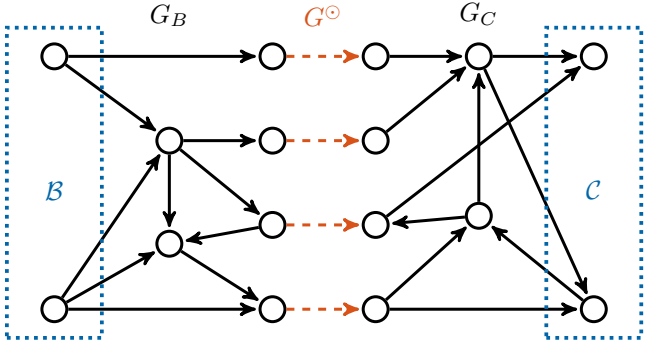

\medskip

\begin{definition} \label{def:separable_net}
	A \emph{separable network} is a network whose matrices have the following block structure:
	\begin{align} \label{eq:def_separable_net}
		G =
        \begin{bmatrix}
            \Gright & G^\odot \\
            0 & \Gleft
        \end{bmatrix},
		&& G^\bullet =
        \begin{bmatrix}
            \Gright & 0 \\
            0 & \Gleft
        \end{bmatrix},
		&& G^\circ =
        \begin{bmatrix}
            0 & G^\odot \\
            0 & 0
        \end{bmatrix}, \nonumber \\
		B =
		\begin{bmatrix}
	        0 & 0 \\
	        0 & \Bleft
	    \end{bmatrix},
		&& C =
        \begin{bmatrix}
            \Cright & 0 \\
            0 & 0
        \end{bmatrix}.
	\end{align}
\end{definition}

\medskip

We have the following important property: on separable networks, global and local identifiability are equivalent.

\medskip

\begin{proposition} \label{prop:separable_global_decoupled}
    A separable network is locally identifiable at $G$ at if and only if it is globally identifiable at $G$.
\end{proposition}

\medskip

\begin{proof}
    Consider a separable network: its matrices have the block structure described in \eqref{eq:def_separable_net}.
    From Definition \ref{def:global_identif}, the network is generically identifiable at $G$ if, for all $\tilde{G}^\circ$ with same zero entries as $G^\circ$, there holds
    \begin{align} \label{eq:global_identif_sep_net}
		C \, (I - \tilde{G})^{-1} \, B = C \, (I- G)^{-1} \,  B \Rightarrow \tilde{G}^\circ = G^\circ
	\end{align}
	where matrices $C,B,G$ and $G^\circ$ have the block structure of \eqref{eq:def_separable_net}, and
		$
        \tilde{G} =
        \begin{bmatrix}
            \Gright & \tilde G^\odot \\
            0 & \Gleft
    	\end{bmatrix}.
    	$
		
		\smallskip
		
    Developing \eqref{eq:global_identif_sep_net} from the block structure of \eqref{eq:def_separable_net} yields
    \begin{align*}
        \Cright \, T(\Gright) \, \tilde G^\odot \, T(\Gleft) \, \Bleft =
        \Cright \, T(\Gright) \, G^\odot \, T(\Gleft) \, \Bleft \\
		\Rightarrow
        \tilde G^\odot = G^\odot,
    \end{align*}
    and bringing out common terms gives
    \begin{align*}
        \Cright  \, T(\Gright) \, \underbrace{(\tilde{G}^\odot - G^\odot)}_{\triangleq \Delta} \, T(\Gleft) \, \Bleft = 0 \\
		\Rightarrow \underbrace{\tilde{G}^\odot - G^\odot}_{= \Delta} = 0,
    \end{align*}
	which is exactly what we obtain by developing \eqref{eq:CTDeltaTB_Delta_net} with the block structure of \eqref{eq:def_separable_net}. To conclude, we know from Theorem \ref{thm:network_local_identif} that \eqref{eq:CTDeltaTB_Delta_net} is a necessary and sufficient condition for local identifiability.
\end{proof}

\medskip

Since $m_\circ = n_B n_C$, $K$ introduced in \eqref{eq:def_K} is a square matrix, hence $\rank K = m_\circ$ is equivalent to $\det K \neq 0$. Theorem \ref{thm:network_local_identif} can then be rewritten in terms of the determinant. In addition, since we work on separable networks, global and local identifiability are equivalent, hence the theorem below characterizes global identifiability.

\medskip

\begin{theorem} \label{thm:identif_separable_detK}
	Consider a separable network.
	Exactly one of the two following holds:
	\begin{enumerate}[(i)]
	    \item $\det K \neq 0$ for almost all $G$ and $G^\circ$ is globally identifiable at almost all $G$;
		\item $\det K = 0$ for all $G$ and $G^\circ$ is globally \emph{non-}identifiable at all $G$.
	\end{enumerate}
	Moreover, $\det K = m_\circ$ is equivalent to the following implication holding for all $\Delta$ with same zero entries as $G^\circ$:
	\begin{align} \label{eq:CTDeltaTB_Delta_net_separable}
		\Cright \, T(\Gright) \, \Delta \, T(\Gleft) \, \Bleft = 0 \Rightarrow \Delta = 0.
	\end{align}
\end{theorem}

\medskip


\medskip

\section{COMBINATORIAL CHARACTERIZATION}




\medskip

We are now going to derive a combinatorial characterization based on a re-expression of the determinant of $K$.
The closed-loop transfer function $T_{ji}$ is expressed in terms of transfer functions $G_{ji}$ in the following way. The analytic matrix $T(G)$ can be expanded in the Taylor series
\begin{align} \label{eq:Taylor_T}
    T(G) = (I-G)^{-1} = I + \sum_{k=1}^\infty G^k,
\end{align}
which converges since the spectral radius of $G$ is strictly smaller than one, see Assumption \ref{assumption:eigenvalues_G}) in Section \ref{subsec:framework}.
A classical result in graph theory is that $[G^k]_{ji}$ is the sum of all paths from $i$ to $j$ of length $k$: 
\begin{align} \label{eq:power_of_G}
    [G^k]_{ji} = \sum_{\substack{\text{all $k$-paths}\\ i \rightarrow j}} \underbrace{G_{j*} \ \dots \ G_{*i} }_{k \text{ terms}},
\end{align}

\noindent where the notation $*$ denotes some node of the network, which can be different for each occurence of $*$.
Combining \eqref{eq:Taylor_T} and \eqref{eq:power_of_G} gives the following lemma, which extends Lemma 2 of \cite{cheng2022identifiability} for networks with cycles.

\medskip

\begin{lemma} \cite{cheng2022identifiability}
    Consider the closed-loop transfer matrix $T(G) = (I-G)^{-1}$. Its entries are given as follows:
    \begin{enumerate}
        \item if there is no path from $i$ to $j$, $T_{ji} = 0$
        \item otherwise, $T_{ji}$ is the (possibly infinite) sum of all the paths from $i$ to $j$:
        \begin{align} \label{eq:Tji}
            T_{ji} = \sum_{\substack{\text{all paths}\\ i \rightarrow j}} G_{j*} \ \dots \ G_{*i}.
        \end{align}
    \end{enumerate}
\end{lemma}

\medskip
In the sequel we refer to the unknown transfer functions $G^\circ_{ji}$ as unknown edges $\alpha$.
We develop $\det K$ as the sum over all possible row-column permutations by Leibniz formula\footnote{$T_{\alpha, \sigma_B(\alpha)}$ is the transfer function between $\sigma_B(\alpha)$, and \emph{start node of} edge $\alpha$, and $T_{\sigma_C(\alpha), \alpha}$ is the one between \emph{end node of} edge $\alpha$ and $\sigma_C(\alpha)$.}:
\begin{align} \label{eq:leibniz_formula}
    \det K 
    = \sum_{\sigma \in S} \sgn(\sigma) \prod_{\alpha \in \mathcal{E}^\circ} T_{\sigma_C(\alpha), \alpha}  \ T_{\alpha, \sigma_B(\alpha)},
\end{align}

\noindent where each row-column permutation corresponds to a bijective assignation $\sigma:\mathcal{E}^\circ \rightarrow \mathcal{B} \times \mathcal{C}$, i.e. between unknown edges and (excitation - measurement) pairs. $S$ denotes the set of all such bijective assignations, $\sigma_B(\alpha)$ is the excitation assigned to edge $\alpha$ by assignation $\sigma$ and $\sigma_C(\alpha)$ is its measurement.
The sign $\sgn(\sigma)$ equals $+1$ if the number of transpositions in assignation $\sigma$ is even, and $-1$ otherwise. A transposition is the swap of two elements, and each $\sigma$ is obtained by combining a certain number of transpositions. 

\medskip

The Leibniz formula \eqref{eq:leibniz_formula} can be further developed by plugging \eqref{eq:Tji} in the equation. The analysis is no longer made on the $T_{ji}$, $\det K$ is now expressed in terms of paths of transfer functions $G_{ji}$\footnote{$G_{\alpha, *}$ stands for the transfer function from a node $*$ to \emph{start node of} edge $\alpha$, and $G_{*, \alpha}$ denotes the one from \emph{end node of} edge $\alpha$ to a node $*$.}:

\begin{align*}
    \det K
    = \sum_{\sigma \in S} \sgn(\sigma) \prod_{\alpha \in \mathcal{E}^\circ} \Biggl[
    &\biggl( \sum_{\substack{\text{all paths}\\ \alpha \rightarrow \sigma_C(\alpha)}} G_{\sigma_C(\alpha),*}  \ \dots \  G_{*,\alpha} \biggr) \\
    \cdot &\biggl( \sum_{\substack{\text{all paths}\\ \sigma_B(\alpha) \rightarrow \alpha}} G_{\alpha,*} \ \dots \  G_{*,\sigma_B(\alpha)} \biggr)  \Biggr]
\end{align*}

This expression can be conveniently re-expressed using the notion of \emph{collections of paths}, that we introduce below.

\medskip

\begin{definition} \label{def:pi}
A collection of paths $\pi$ is a set of $m^\circ$ paths $\pi_i$ (one $\pi_i$ for each unknown edge), where each path $\pi_i$ is a sequence of connected edges that starts at an excited node (i.e. in $\mathcal B$), and ends at a measured node (i.e. in $\mathcal C$).
Since we work on separable networks, each $\pi_i$ goes through one and only one unknown edge.

\noindent
We say that $\pi$ is \emph{bijective} if no two paths of $\pi$ start at the same excitation and end at the same measurement.
\end{definition}

\medskip

\noindent Distributing the product over all unknown edges $\alpha$ gives
\begin{align*}
    \det K
    = \sum_{\sigma \in S} \sgn(\sigma) \sum_{\pi \in \Pi_\sigma} \mu(\pi),
\end{align*}
where $\pi$ is a bijective collection of paths, as defined in Definition \ref{def:pi}, $\Pi_\sigma$ is the set of all $\pi$ corresponding to assignation $\sigma$ and $\mu(\pi)$ is the \emph{monomial} corresponding to $\pi$:
\begin{align} \label{eq:detK_sum_sigma_pi}
    \mu(\pi) \triangleq \prod_{\alpha \in \mathcal{E}^\circ} G_{\pi_C(\alpha),*} \ ... \ G_{*,\alpha} \cdot G_{\alpha,*}
     \ ... \ G_{*,\pi_B(\alpha)}
\end{align}
where $\pi_B(\alpha)$ and $\pi_C(\alpha)$ are respectively the starting and ending node of the path going through $\alpha$ in $\pi$, and the transfer functions $G_{ji}$ may appear several times in $\mu(\pi)$, which would mean that it is taken several times by $\pi$.

\medskip

From \eqref{eq:detK_sum_sigma_pi}, we group by monomials of same $\pi$: instead of summing over $\sigma$ and then over $\pi$, we sum directly over $\pi$:
\begin{align} \label{eq:detK_sum_pi}
    \det K
    = \sum_{\pi \in \Pi} \sgn(\pi) \, \mu(\pi),
\end{align}
where $\Pi$ is the set of all bijective collection of paths $\pi$. As $\pi$ derives from an assignation $\sigma$, its monomial has a sign, denoted by $\sgn(\pi) $.


\medskip


Observe that from a same group of edges it may be possible to build different bijective $\pi$ (this will be the case e.g. when two paths cross multiple times), so that these different $\pi$ will have the same monomial $\mu$. Hence we regroup the same monomials together:
\begin{align} \label{eq:detK_sum_mu_noncompact}
    \det K
    = \sum_{\mu \in M} \biggl( \sum_{\pi \in \Pi_\mu} \sgn(\pi) \biggr) \ \mu,
\end{align}
where $M$ is the set of all monomials $\mu$ corresponding to bijective collections of paths, and $\Pi_\mu$ is the set of all bijective collections of paths $\pi$ corresponding to monomial $\mu$. We define the \emph{repetition} $r(\mu) \triangleq \sum_{\pi \in \Pi_\mu} \sgn(\pi)$,
which allows to rewrite \eqref{eq:detK_sum_mu_noncompact} in a compact way:
\begin{align} \label{eq:detK_sum_mu}
    \det K
    = \sum_{\mu \in M} r(\mu) \, \mu.
\end{align}

As seen in Theorem \ref{thm:identif_separable_detK}, $\det K$ encodes the global (non)-identifiability of separable networks. Hence equation \eqref{eq:detK_sum_mu} allows providing a necessary and sufficient combinatorial characterization of global identifiability, in terms of the repetition $r(\mu)$.

\toadd{Provide an example where $r(\mu)$ is easy to compute and understand, and one where it is complicated.}

\medskip

\begin{theorem} \label{thm:combinatorial_charac}
    Consider a separable network. It is generically globally identifiable if and only if there is at least one monomial $\mu \in M$ such that its repetition $r(\mu) \neq 0$.
\end{theorem}

\medskip

\begin{proof}
    From Theorem \ref{thm:identif_separable_detK}, we know that a separable network is generically globally identifiable if and only if $\det K \neq 0$ for almost all $G$.
    Equation \eqref{eq:detK_sum_mu} expresses $\det K$ as a sum over all monomials $\mu$, weighted by their repetition $r(\mu)$.
    A sum of distinct monomials is generically nonzero if and only if at least one of the monomials has a nonzero coefficient.
    Therefore, for this sum to be nonzero for almost all $G$, at least one $\mu$ must have a nonzero repetition $r(\mu)$
\end{proof}

\vspace{-1mm}

Since the decoupled network is a particular case of separable network, this condition can be formulated on the decoupled network of a network $G$ of any topology, \emph{not only separable networks}. And generic decoupled identifiability is necessary for generic global identifiability \cite{legat2021path}, hence the necessary condition of Theorem \ref{thm:combinatorial_charac}, formulated on decoupled network $\hat G$, applies to a network $G$ of any topology.

\medskip

Theorem \ref{thm:combinatorial_charac} gives a necessary and sufficient combinatorial characterization, but building an efficient algorithm to check this condition 
remains an open question.

\medskip

\section{CONCLUSION}
\medskip

This work was motivated by one main open question: determining path-based conditions for generic identifiability of networked systems.

\medskip

Introducing separable networks allowed to address global identifiability of such graphs using the algebraic tools of local identifiability. A new approach led to a necessary and sufficient combinatorial characterization of identifiability, in terms of existence of paths and conditions on their parity.



\medskip

Furthermore, this necessary condition not only applies to separable networks, but to networks of \emph{any topology.} It follows from the fact that the decoupled network is a particular case of separable networks, and generic decoupled identifiability is necessary for generic global identifiability.

\medskip

A further open question is whether our necessary and sufficient combinatorial characterization can be algorithmically checked.
Also, we would like to establish the equivalence (or not) between the notions introduced here and local identifiability for a general network.




\smallskip

\section{ACKNOWLEDGMENTS}

\smallskip

The authors gratefully acknowledge Federica Garin (Inria, Gipsa-lab, France) for the interesting, insightful discussions.

\smallskip

\bibliographystyle{ieeetr}
\bibliography{0-cdc23}

\begin{thebibliography}{10}

\bibitem{weerts2015identifiability}
H.~H. Weerts, A.~G. Dankers, and P.~M. Van~den Hof, ``Identifiability in
  dynamic network identification,'' {\em IFAC-PapersOnLine}, 2015.

\bibitem{weerts2018single}
H.~Weerts, P.~M. Van~den Hof, and A.~Dankers, ``Single module identifiability
  in linear dynamic networks,'' in {\em 2018 IEEE Conference on Decision and
  Control (CDC)}, pp.~4725--4730, IEEE, 2018.

\bibitem{gevers2018practical}
M.~Gevers, A.~S. Bazanella, and G.~V. da~Silva, ``A practical method for the
  consistent identification of a module in a dynamical network,'' {\em
  IFAC-PapersOnLine}, vol.~51, no.~15, pp.~862--867, 2018.

\bibitem{weerts2018prediction}
H.~H. Weerts, P.~M. Van~den Hof, and A.~G. Dankers, ``Prediction error
  identification of linear dynamic networks with rank-reduced noise,'' {\em
  Automatica}, vol.~98, pp.~256--268, 2018.

\bibitem{gevers2018identifiability}
M.~Gevers, A.~S. Bazanella, and G.~A. Pimentel, ``Identifiability of dynamical
  networks with singular noise spectra,'' {\em IEEE Transactions on Automatic
  Control}, vol.~64, no.~6, pp.~2473--2479, 2018.

\bibitem{van2019local}
P.~M. Van~den Hof, K.~R. Ramaswamy, A.~G. Dankers, and G.~Bottegal, ``Local
  module identification in dynamic networks with correlated noise: the full
  input case,'' in {\em 2019 IEEE 58th Conference on Decision and Control
  (CDC)}, pp.~5494--5499, IEEE, 2019.

\bibitem{ramaswamy2020local}
K.~R. Ramaswamy and P.~M. Vandenhof, ``A local direct method for module
  identification in dynamic networks with correlated noise,'' {\em IEEE
  Transactions on Automatic Control}, 2020.

\bibitem{bazanella2017identifiability}
A.~S. Bazanella, M.~Gevers, J.~M. Hendrickx, and A.~Parraga, ``Identifiability
  of dynamical networks: which nodes need be measured?,'' in {\em IEEE 56th
  Annual Conference on Decision and Control (CDC)}, 2017.

\bibitem{weerts2018identifiability}
H.~H. Weerts, P.~M. Van~den Hof, and A.~G. Dankers, ``Identifiability of linear
  dynamic networks,'' {\em Automatica}, vol.~89, pp.~247--258, 2018.

\bibitem{hendrickx2018identifiability}
J.~M. Hendrickx, M.~Gevers, and A.~S. Bazanella, ``Identifiability of dynamical
  networks with partial node measurements,'' {\em IEEE Transactions on
  Automatic Control}, vol.~64, no.~6, pp.~2240--2253, 2018.

\bibitem{shi2020excitation}
S.~Shi, X.~Cheng, and P.~M. Van~den Hof, ``Excitation allocation for generic
  identifiability of a single module in dynamic networks: A graphic approach,''
  in {\em Proc. 21st IFAC World Congress}, 2020.

\bibitem{cheng2019allocation}
X.~Cheng, S.~Shi, and P.~M. Van~den Hof, ``Allocation of excitation signals for
  generic identifiability of dynamic networks,'' in {\em Proceedings of the
  IEEE Conference on Decision and Control}, 2019.

\bibitem{cheng2021allocation}
X.~Cheng, S.~Shi, and P.~M. Van~den Hof, ``Allocation of excitation signals for
  generic identifiability of linear dynamic networks,'' {\em IEEE Transactions
  on Automatic Control}, vol.~67, no.~2, pp.~692--705, 2021.

\bibitem{dreef2022excitation}
H.~Dreef, S.~Shi, X.~Cheng, M.~Donkers, and P.~M. Van~den Hof, ``Excitation
  allocation for generic identifiability of linear dynamic networks with fixed
  modules,'' {\em IEEE Control Systems Letters}, vol.~6, pp.~2587--2592, 2022.

\bibitem{shi2020generic}
S.~Shi, X.~Cheng, and P.~M. Van~den Hof, ``Generic identifiability of
  subnetworks in a linear dynamic network: the full measurement case,'' {\em
  arXiv preprint arXiv:2008.01495}, 2020.

\bibitem{van2018topological}
H.~J. van Waarde, P.~Tesi, and M.~K. Camlibel, ``Topological conditions for
  identifiability of dynamical networks with partial node measurements,'' {\em
  IFAC-PapersOnLine}, vol.~51, no.~23, pp.~319--324, 2018.

\bibitem{van2019necessary}
H.~J. Van~Waarde, P.~Tesi, and M.~K. Camlibel, ``Necessary and sufficient
  topological conditions for identifiability of dynamical networks,'' {\em IEEE
  Transactions on Automatic Control}, 2019.

\bibitem{bazanella2019network}
A.~S. Bazanella, M.~Gevers, and J.~M. Hendrickx, ``Network identification with
  partial excitation and measurement,'' in {\em 58th Conference on Decision and
  Control (CDC)}, pp.~5500--5506, IEEE, 2019.

\bibitem{cheng2022identifiability}
X.~Cheng, S.~Shi, I.~Lestas, and P.~M. Van~den Hof, ``Identifiability in
  dynamic acyclic networks with partial excitations and measurements,'' {\em
  arXiv preprint arXiv:2201.07548}, 2022.

\bibitem{shi2021single}
S.~Shi, X.~Cheng, and P.~M. Van~den Hof, ``Single module identifiability in
  linear dynamic networks with partial excitation and measurement,'' {\em IEEE
  Transactions on Automatic Control}, 2021.

\bibitem{cheng2021necessary}
X.~Cheng, S.~Shi, I.~Lestas, and P.~M. Van~den Hof, ``A necessary condition for
  network identifiability with partial excitation and measurement,'' {\em arXiv
  preprint arXiv:2105.03187}, 2021.

\bibitem{mapurunga2022identifiability}
E.~Mapurunga, M.~Gevers, and A.~S. Bazanella, ``Identifiability of dynamic
  networks: the essential r$\backslash$\^{} ole of dources and dinks,'' {\em
  arXiv preprint arXiv:2211.10825}, 2022.

\bibitem{mapurunga2022excitation}
E.~Mapurunga, M.~Gevers, and A.~S. Bazanella, ``Excitation and measurement
  patterns for the identifiability of directed acyclic graphs,'' in {\em 2022
  IEEE 61st Conference on Decision and Control (CDC)}, pp.~1616--1621, IEEE,
  2022.

\bibitem{legat2020local}
A.~Legat and J.~M. Hendrickx, ``Local network identifiability with partial
  excitation and measurement,'' in {\em 2020 59th IEEE Conference on Decision
  and Control (CDC)}, pp.~4342--4347, IEEE, 2020.

\bibitem{legat2021path}
A.~Legat and J.~M. Hendrickx, ``Path-based conditions for local network
  identifiability -- full version,'' {\em arXiv preprint}, 2021.

\bibitem{van2019topology_heterogeneous}
H.~J. van Waarde, P.~Tesi, and M.~K. Camlibel, ``Topology identification of
  heterogeneous networks of linear systems,'' in {\em 58th Conference on
  Decision and Control (CDC)}, pp.~5513--5518, IEEE, 2019.

\bibitem{van2019topology_lyapunov}
H.~J. van Waarde, P.~Tesi, and M.~K. Camlibel, ``Topology reconstruction of
  dynamical networks via constrained lyapunov equations,'' {\em IEEE
  Transactions on Automatic Control}, vol.~64, pp.~4300--4306, 2019.

\bibitem{kivits2022local}
E.~L. Kivits and P.~M. Van~den Hof, ``Local identification in diffusively
  coupled linear networks,'' in {\em 2022 IEEE 61st Conference on Decision and
  Control (CDC)}, pp.~874--879, IEEE, 2022.

\bibitem{davison1977connectability}
E.~J. Davison, ``Connectability and structural controllability of composite
  systems,'' {\em Automatica}, vol.~13, no.~2, pp.~109--123, 1977.

\bibitem{dion2003generic}
J.-M. Dion, C.~Commault, and J.~Van~der Woude, ``Generic properties and control
  of linear structured systems: a survey,'' {\em Automatica}, 2003.

\bibitem{matlab}
A.~Legat and J.~M. Hendrickx, {\em Identifiability test}.
\newblock \url{https://github.com/alegat/identifiable}.

\end{thebibliography}


\end{document}